\documentclass [12pt]{article}
\title {Supersimplicity and countable reducts of a unidimensional hypersimple theory}
\author {Ziv Shami}
\newtheorem {theorem}{Theorem}[section]
\newtheorem {lemma}[theorem]{Lemma}
\newtheorem {definition}[theorem]{Definition}

\newtheorem {fact}[theorem]{Fact}

\newtheorem {remark}[theorem]{Remark}

\newtheorem {proposition}[theorem]{Proposition}

\newtheorem {notation}[theorem]{Notation}
\newtheorem {claim}[theorem]{Claim}

\def\proof {\noindent \textbf{Proof:} }

\def\qed {$\ \ \ \ \Box$}

\def\Ceq  {\CC^{eq}}

\newsavebox{\indbin}
\savebox{\indbin}{\begin{picture}(0,0)
\newlength{\gnu}
\settowidth{\gnu}{$\smile$} \setlength{\unitlength}{.5\gnu} \put(-1,-.65){$\smile$}
\put(-.25,.1){$|$}
\end{picture}}
\newcommand{\nonfork}[3]
{\mbox{$\begin{array}{ccc} \mbox{$#1$} & \usebox{\indbin} & \mbox{$#2$} \\
        & \mbox{$#3$} &
\end{array}$}}
\newcommand{\nonforkempty}[2]
{\mbox{$\begin{array}{ccc} \mbox{$#1$} & \usebox{\indbin} & \mbox{$#2$}
\end{array}$}}
\newcommand{\fork}[3]
{\mbox{$\begin{array}{ccc} \mbox{$#1$} & \!\mbox{$\!\!\not\!\:\usebox{\indbin}$} & \mbox{$#2$} \\
        & \mbox{$#3$} &
\end{array}$}}

\newsavebox{\sindbin}
\savebox{\sindbin}{\begin{picture}(0,0)
\newlength{\sgnu}
\settowidth{\sgnu}{$\smile$} \setlength{\unitlength}{.5\sgnu} \put(-1,-.65){$\smile$}
\put(-.25,.1){$|s$}
\end{picture}}
\newcommand{\snonfork}[3]
{\mbox{$\begin{array}{ccc} \mbox{$#1$} & \usebox{\sindbin} & \mbox{$#2$} \\
        & \mbox{$#3$} &
\end{array}$}}

\newcommand{\sfork}[3]
{\mbox{$\begin{array}{ccc} \mbox{$#1$} & \!\mbox{$\!\!\not\!\:\usebox{\sindbin}$} & \mbox{$#2$} \\
        & \mbox{$#3$} &
\end{array}$}}

\newsavebox{\starindbin}
\savebox{\starindbin}{\begin{picture}(0,0)
\newlength{\stargnu}
\settowidth{\stargnu}{$\smile$} \setlength{\unitlength}{.5\stargnu} \put(-1,-.65){$\smile$}
\put(-.25,.1){$|*$}
\end{picture}}

\newsavebox{\rindbin}
\savebox{\rindbin}{\begin{picture}(0,0)
\newlength{\rgnu}
\settowidth{\rgnu}{$\smile$} \setlength{\unitlength}{.5\rgnu} \put(-1,-.65){$\smile$}
\put(-.25,.1){$|-$}
\end{picture}}

\newsavebox{\qindbin}
\savebox{\qindbin}{\begin{picture}(0,0)
\newlength{\qgnu}
\settowidth{\qgnu}{$\smile$} \setlength{\unitlength}{.5\qgnu} \put(-1,-.65){$\smile$}
\put(-.25,.1){$|_{qf}$}
\end{picture}}

\def\card #1 {{\vert #1 \vert}}

\def\CC {{\cal C}}
\def\RR {{\cal R}}

\def\PP {{\cal P}}
\def\UU {{\cal U}}

\def\proofof #1 {\noindent\textbf{Proof of Theorem \ref{#1}}}

\usepackage{amsfonts}

\begin{document}
\maketitle

\begin{abstract}
We show that a hypersimple unidimensional theory that has a club of reducts, in the partial order of all countable reducts, that are
coordinatized in finite rank, is supersimple.
\end{abstract}


\section{Introduction} In this paper we suggest an approach to the problem on
supersimplicity of unidimensional hypersimple theories. The problem has been answered in the
affirmative in the following cases. In [H], for any stable theory, in [S1] for any countable theory (this improved an earlier result for the case of a countable theory with the wnfcp [P]) and
in [S2] it has been proved for any (possibly uncountable) non-s-essentially 1-based  theory (roughly, a theory that is far from being 1-based).

It is easy to see that supersimplicity of a theory is determined by (the supersimplicity of) the family of its countable reducts.
Therefore, it is natural to try and reflect properties of the given unidimensional hypersimple theory to countable reducts. Clearly, unidimensionality is
not preserved under reducts. On the other hand, easily any unidimensional hypersimple theory is coordinatized in finite rank (see Definition 5.1).
In this paper we show that supersimplicity of any (possibly uncountable) unidimensional hypersimple theory follows from coordinatization in finite rank
of sufficiently many countable reducts of it.

We thank Ehud Hrushovski for discussions on this topic and for allowing us to include a remark of him about elimination of hyperimaginaries in reducts (section 3).
We will assume basic knowledge of simple theories as in [K1],[KP],[HKP]. A good text book on simple theories that covers much more is [W]. The notations
are standard, and throughout the paper we work in a large saturated model $\CC$ of a complete first-order theory $T$ in a language $L$.

\section{Preliminaries} In this section $T$ is assumed to be simple. We quote several known facts
that we will apply.

\subsection{Almost internality, analyzability and unidimensionality}
In this subsection we work in $\CC$ with hyperimaginaries unless otherwise stated; if $T$ is hypersimple (i.e. simple and eliminates hyperimaginaries) and we work in $\Ceq$ we get equivalent definitions.\\

\noindent In this subsection, $\PP$ denotes an $A$-invariant family of partial types and $p$ a partial type over $A$. We say that \em $p$ is (almost-)
$\PP$-internal \em for every realization $a$ of $p$ there exists $b$ with
$\nonfork{a}{b}{A}$ such that for some tuple $c$ of realizations of partial types in $\PP$ over
$Ab$ we have $a \in dcl(b,c)$ (respectively, $a \in acl(b,c)$). We say that \em $p$ is analyzable
in $\PP$ \em if for any $a\models p$ there exists a sequence $I=\langle a_i \vert i\leq \alpha \rangle$
such that $a_\alpha=a$ and $tp(a_i/\{a_j \vert j<i\}\cup
A)$ is almost $\PP$-internal for every $i\leq \alpha$.\\

\noindent First, the following fact is straightforward.

\begin{fact}\label{some internal}
1) Assume $tp(a_i/A)$ are (almost) $\PP$-internal for $i<\alpha$. Then $tp(\langle a_i \vert
i<\alpha \rangle/A)$ is (respectively, almost) $\PP$-internal. Thus, if $tp(a_i/A)$ are analyzable in $\PP$ for
$i<\alpha$. Then $tp(\langle a_i \vert i<\alpha \rangle/A)$ is analyzable in $\PP$.

\noindent 2) If $tp(a/A)$ almost $\PP$-internal, so is $tp(a/B)$ for any set $B\supseteq A$.
\end{fact}

\noindent $T$ is said to be \em unidimensional \em if whenever $p$ and $q$ are complete
non-algebraic types, $p$,$q$ are non-orthogonal.\\

We will also need the following easy Fact.

\begin{fact}\label{supersimple definable}
Work in $\CC$ (without hyperimaginaries). Let $p\in S(\emptyset)$ and let $\theta\in L$. Assume $p$ is analyzable in $\theta$. Then $p$
is analyzable in $\theta$ in finitely many steps. In particular, if $T$ is a hypersimple
unidimensional theory and there exists a non-algebraic supersimple definable set, then $T$ has
finite $SU$-rank, i.e. every complete type has finite $SU$-rank. In fact, for every given sort
there is a finite bound on the $SU$-rank of all types in that sort, equivalently the global
$D$-rank of any sort is finite.
\end{fact}

Another useful fact is the following.

\begin{fact}\label{elimination_exists_infty} $[S1]$
Let $T$ be any unidimensional simple theory. Then $T$ eliminates $\exists^\infty$.
\end{fact}

\subsection{The forking topology, EPFO and PCFT}
The forking topology is introduced in [S0] and is a variant of Hrushovski's and Pillay's topologies
from [H0] and [P], respectively. In this section $T$ is assumed to be simple and we work in $\CC$.

\begin{definition} \label{tau definition}\em
Let $A\subseteq \CC$ and let $x$ be a finite tuple of variables.\\
1) An invariant set $\UU$ over $A$ is said to be \em a basic $\tau^f$-open set over $A$ \em if
there is $\phi(x,y)\in L(A)$ such that $$\UU=\{a \vert \phi(a,y)\ \mbox{forks\ over}\ A \}.$$ Note
that the family of basic $\tau^f$-open sets over $A$ is closed under finite intersections, thus
forms a basis for a unique topology on $S_x(A)$. An open set in this topology is called a
$\tau^f$-open set over $A$ or a forking-open set over $A$.

\noindent 2) An invariant set $\UU$ over $A$ is said to be \em a basic $\tau^f_\infty$-open set
over $A$ \em if $\UU$ is a type-definable $\tau^f$-open set over $A$. The family of basic
$\tau^f_\infty$-open sets over $A$ is a basis for a unique topology on $S_x(A)$. An open set in
this topology is called a $\tau^f_\infty$-open set over $A$.
\end{definition}

\begin{remark}\label{exists generic}\em
The $\tau^f_\infty$-topology and in particular the $\tau^f$-topology on $S_x(A)$ refines the Stone-topology of $S_x(A)$ for all
$x,A$.
\end{remark}

We will apply the following Fact.

\begin{fact}\label{tau extensions}
Let $\UU$ be a $\tau^f$-open set over $\emptyset$ and let $A$ be any set. Then $\UU$ is
$\tau^f$-open over $A$.
\end{fact}

Recall the following definition from [S0].

\begin{definition}\label{projection closed}\em
We say that \em the $\tau^f$-topologies over $A$ are closed under projections ($T$ is PCFT over
$A$) \em if for every $\tau^f$-open set $\UU(x,y)$ over $A$ the set $\exists y \UU(x,y)$ is a
$\tau^f$-open set over $A$. We say that \em the $\tau^f$-topologies are closed under projections
($T$ is PCFT) \em if they are such over every set $A$.
\end{definition}

In [BPV, Proposition\ 4.5] the authors proved the following equivalence which, for convenience, we
will use as a definition (their definition involves extension with respect to pairs of models of
$T$).

\begin {definition}\label {foext}\em
We say that the extension property is first-order in $T$ or $T$ is EPFO iff for every formulas
$\phi(x,y),\psi(y,z)\in L$ the relation $Q_{\phi,\psi}$ defined by: $$Q_{\phi,\psi}(a)\mbox{\ iff}\
\phi(x,b)\mbox{ doesn't\ fork\ over}\ a\ \mbox{for\ every}\ b\models \psi(y,a)$$ is type-definable
(here $a$ can be an infinite tuple from $\CC$ whose sorts are fixed).
\end {definition}

\begin{fact}$[S1, Corollary\ 3.13]$\label{ext pcft}
Suppose the extension property is first-order in $T$. Then $T$ is PCFT.
\end{fact}

We say that an $A$-invariant set $\UU$ \em has finite $SU$-rank \em if $SU(a/A)<\omega$ for all
$a\in\UU$, and \em has bounded finite $SU$-rank \em if there exists $n<\omega$ such that
$SU(a/A)\leq n$ for all $a\in\UU$. The existence of a $\tau^f$-open set of bounded finite $SU$-rank
implies the existence of an $SU$-rank 1 formula (i.e. a weakly-minimal formula):

\begin{fact}\label{tau bounded SU}$[S0, Proposition\ 2.13]$
Let $\UU$ be an unbounded $\tau^f$-open set over some set $A$. Assume $\UU$ has bounded finite
$SU$-rank. Then there exist a set $B\supseteq A$ with $\vert B\backslash A\vert<\omega$ and
$\theta(x)\in L(B)$ of $SU$-rank 1 such that $\theta^\CC\subseteq \UU\cup acl(B)$.
\end{fact}

Now, recall the following two facts and their corollary. First, let $\PP^{SU\leq 1}$ denote the
class of complete real types over sets of size $\leq \vert T \vert$, of $SU$-rank $\leq 1$.

\begin{fact} \label {fact 1} $[P1]$
Let $T$ be a simple theory that eliminates $\exists^\infty$. Moreover, assume every
type is analyzable in $\PP^{SU\leq 1}$. Then the extension property is first-order in $T$.
\end{fact}
\noindent For a more general statement, see [S1, Lemma 3.7].

\subsection{Stable independence and stable SU-rank}
In this subsection $T$ is assumed to be simple and we work in $\CC$.\\

First we recall the notion of stable independence.

\begin{definition}\label{stable_dep}\em
Let $a\in \CC$, $A\subseteq B\subseteq \CC$. We say that \em $a$ is stably-independent from $B$
over $A$ \em if for every stable $\phi(x,y)\in L$, if $\phi(x,b)$ is over $B$ and $a'\models
\phi(x,b)$ for some $a'\in dcl(Aa)$, then $\phi(x,b)$ doesn't divide over $A$. In this case we
denote it by $\snonfork{a}{B}{A}$.
\end{definition}

The notion of stable $SU$-rank is defined via stable dependence.

\begin{definition}\em
1) For $a\in \CC$ and $A\subseteq \CC$ the $SU_s$-rank is defined by induction on $\alpha$: if
$\alpha=\beta+1$, then $SU_s(a/A)\geq\alpha$ if there exists $B\supseteq A$ such that
$\sfork{a}{B}{A}$ and $SU_s(a/B)\geq\beta$. For limit $\alpha$, $SU_s(a/A)\geq\alpha$ if
$SU_s(a/A)\geq\beta$ for all $\beta<\alpha$.

\noindent 2) Let $\UU$ be an $A$-invariant set. We write $SU_s(\UU)=\alpha$ (the $SU_s$-rank of
$\UU$ is $\alpha$) if $Max\{SU_s(p) \vert p\in S(A), p^\CC\subseteq\UU\}=\alpha$. We say that \em
$\UU$ has bounded finite $SU_s$-rank \em if for some $n<\omega$, $SU_s(\UU)=n$. Note that the
$SU_s$-rank of $\UU$ might, a priori, depend on the choice of the set $A$ over which $\UU$ is
invariant.
\end{definition}

The following rank is a variation of stable $SU$-rank; it is non-increasing in extensions.

\begin{definition}\em
1) For $a\in \CC$ and $A\subseteq \CC$ the $SU_{se}$-rank is defined by induction on $\alpha$: if
$\alpha=\beta+1$, $SU_{se}(a/A)\geq \alpha$ if there exist $B_1\supseteq B_0\supseteq A$ such that
$\sfork{a}{B_1}{B_0}$ and $SU_{se}(a/B_1)\geq\beta$. For limit $\alpha$, $SU_{se}(a/A)\geq\alpha$
if $SU_{se}(a/A)\geq\beta$ for all $\beta<\alpha$.

\noindent 2) Let $\UU$ be an $A$-invariant set. We write $SU_{se}(\UU)=\alpha$ (the $SU_{se}$-rank
of $\UU$ is $\alpha$) if $Max\{SU_{se}(p) \vert p\in S(A), p^\CC\subseteq\UU\}=\alpha$. We say that
\em $\UU$ has bounded finite $SU_{se}$-rank \em if for some $n<\omega$, $SU_{se}(\UU)=n$.
\end{definition}

\begin{remark} \label {SU_s SU_se}\em
Note that $SU_{se}(a/B)\leq SU_{se}(a/A)$ for all $a\in \CC$ and $A\subseteq B\subseteq\CC$ (this
is the reason for introducing $SU_{se})$. Also, clearly $SU_s(a/A)\leq SU_{se}(a/A)\leq SU(a/A)$
for all $a,A$. Clearly $SU_{se}(a/A)=0$ iff $SU_s(a/A)=0$ iff $a\in acl(A)$ for all $a,A$.
\end{remark}

We will apply the following easy fact.

\begin{fact}\label{stable forking}
For $a\in\CC$ and $A\subseteq B\subseteq \CC$, assume $tp(a/B)$ doesn't fork over $acl(aA)\cap
acl(B)$ and $\fork{a}{B}{A}$. Then $\sfork{a}{B}{A}$.
\end{fact}

\section{Elimination of hyperimaginaries in reducts}
In this section we include a remark by Ehud Hrushovski that allowed us to remove the assumption that the reducts eliminate hyperimaginaries (in the main theorem).
Here $T$ denotes any complete theory in a language $L$  and we work in $\CC$.

\begin{definition}
A reduct $T^-$ of $T$ to a sublanguage $L^-\subseteq L$ is said to be \em  $E$-closed \em if for
every $L^-$-definable sets $D_1\vdash D_2$ on $S^2$ (for some sort $S$ of $L^-$) there exists a
definable equivalence relation $E^-\in L^-$ satisfying $D_1\vdash E^-\vdash D_2$, provided that
there exists such definable equivalence relation in $L$.
\end{definition}

\noindent For a partial order $(P,\leq)$, a subset $A\subseteq P$ is called a \em club \em in $(P,\leq)$, if
$A$ unbounded in $(P,\leq)$, that is,  above any element of $P$ there is an element of
$A$, and $A$ is closed in $(P,\leq)$, that is,  for any chain $C\subseteq A$, if $a\in P$ is the supremum of $C$
 (i.e. $a$ is an upper bound of $C$ and $a$ is smaller then any other upper bound of $C$)
then $a\in A$.

\begin{notation}
Let  $T^-$ be a reduct of $T$ to $L^-$. The size of the reduct $T^-$ is just $\vert T^-\vert$. Let $\lambda$ be any infinite cardinal (or $\infty$).
Let $(\RR^\lambda_{T},\leq_{T})$ be the partial order of all reducts of $T$ of size $\le\lambda$, where the order is just inclusion
(of the sublanguages of the reducts, i.e. of both the set of sorts and the set of formulas).  It will be convenient to consider the isomorphic
partial order  $(\RR^\lambda_{\CC},\leq_{\CC})$ of all the (saturated) model reducts of $\CC$ to a sublanguage of $L$ size  $\le\lambda$.
\end{notation}

\begin{claim}\label{reduct_elim}
Let $T$ be any complete $L$-theory that eliminates hyperimaginaries.
\noindent 1) Let $T^-$ be an $E$-closed reduct of $T$. Then $T^-$ eliminates hyperimaginaries.
\noindent 2) The set of $E$-closed reducts of $T$ is a club in  $(\RR^\infty_{T},\leq_{T})$.
Given any infinite $\lambda\leq \vert L\vert$, the set of $E$-closed reducts of $T$ of size
$\leq\lambda$ is a club in $(\RR^\lambda_{T},\leq_{T})$.
\end{claim}

\proof 1) Say $T^-$ is the reduct of $T$ to $L^-$ and so $\CC\vert L^-$ is a saturated model of $T^-$. We claim that the hyperimaginaries of $T^-$ are
eliminated, namely: for every type-definable equivalence relation $E^-$ of $T^-$ on a complete type $p^-$ of $T^-$ over $\emptyset$,
 there are definable equivalence relations $E_i^-\in L^-$ such that $E^-$ is
equivalent to $\bigwedge_i E_i^-$ on $p^-$. Indeed, let $E^-=E^-(x,x')$ and $p^-=p^-(x)$ be such. Let $\phi^-_i(x,x')\in L^-$ be such that
$E^-(x,x')=\bigwedge_i\phi^-_i(x,x')$. Let $p$ be any complete type of $T$ over $\emptyset$ that extends $p^-$.
By elimination of hyperimaginaries in $T$, there are $E_j(x,x')\in L$ such that $\bigwedge_j E_j(x,x')$ is equivalent to $E^-(x,x')$ on $p^\CC$.
Now, by compactness, for any $i$ there is $j(i)$ such that $E_{j(i)}(x,x')\vdash \phi^-_i(x,x')$ on $p^\CC$, likewise, for every $j$ there exists
$k(j)$ such that  $\phi^-_{k(j)}(x,x')\vdash E_{j}(x,x')$ on $p^\CC$. As $T^-$ is an $E$-closed reduct of $T$, for every $i$,
there is a definable equivalence relation $E^-_i\in L^-$ such that $\phi^-_{k(j(i))}(x,x')\vdash E^-_i(x,x')\vdash \phi^-_i(x,x')$ on $p^\CC$ (using compactness).
We conclude that $E^-$ is equivalent to $\bigwedge_i E^-_i(x,x')$ on $p^\CC$ and thus on $p^-\CC$ as well (as $E^-$ and $E^-_i$ are all invariant under automorphisms of $\CC\vert L^-$).
 2) is immediate.

\section{Dichotomies for $\emptyset$-invariant families of rank 1 types}

Here we verify the following extension of [S2, Corollary 2.13] to a general
$\emptyset$-invariant family of $SU$-rank 1 types. In this section $T$ is assumed to be a simple theory with elimination of imaginaries.\\

We first recall some basic definitions from [S1].

\begin{definition}\em
A family $$\Upsilon=\{\Upsilon_{x,A} \vert\ x \mbox{ is a finite sequence of variables and }
A\subset \CC \mbox{ is small}\}$$ is said to be \em a projection closed family of topologies \em if
each $\Upsilon_{x,A}$ is a topology on $S_x(A)$ that refines the Stone-topology on $S_x(A)$, this
family is invariant under automorphisms of $\CC$ and change of variables by variables of the same
sort, the family is closed under product by the full Stone spaces $S_y(A)$ (where $y$ is a disjoint
tuple of variables) and closed by projections, namely whenever $\UU(x,y)\in \Upsilon_{xy,A}$,
$\exists y\UU(x,y)\in\Upsilon_{x,A}$.
\end{definition}

We will be interested in the case  $\Upsilon= \tau^f$, where $T$ is a PCFT theory.
From now on fix a general projection closed family $\Upsilon$ of topologies.

\begin{definition}\label {def ess-1-based}\em
1) A type $p\in S(A)$ is said to be \em s-essentially 1-based over $A_0\subseteq A$ (essentially 1-based over $A_0\subseteq A$) by means of
$\Upsilon$ \em if for every finite tuple $\bar c$ from $p$ and for every (respectively, type-definable)
$\Upsilon$-open set $\UU$ over $A\bar c$, with the property that $a$ is independent from $A$ over
$A_0$ for every $a\in \UU$, the set $\{a\in \UU \vert\ Cb(a/A\bar c)\not\in bdd(aA_0)\}$ is nowhere
dense in the Stone-topology of $\UU$. We say $p\in S(A)$ is \em s-essentially 1-based (essentially 1-based) by means of
$\Upsilon$ \em if $p$ is s-essentially 1-based (respectively, essentially 1-based) over $A$ by means of $\Upsilon$.\\
2) Let $V$ be an $A_0$-invariant set and let $p\in S(A_0)$. We say that $p$ is \em analyzable in
$V$ by s-essentially 1-based (by essentially 1-based) types by means of $\Upsilon$ \em if there exists $a\models p$ and there
exists a sequence $(a_i\vert\ i\leq\alpha)\subseteq dcl(A_0a)$ with $a_\alpha=a$ such that
$tp(a_i/A_0\cup\{a_j\vert j<i\})$ is $V$-internal and s-essentially 1-based (respectively, essentially 1-based) over $A_0$ by means of
$\Upsilon$ for all $i\leq\alpha$.
\end{definition}

\begin{theorem}\label{empty_inv_su1_dich_pcft}
Let $T$ be any countable hypersimple theory with PCFT. Let $\PP_0$ be an $\emptyset$-invariant
family of $SU$-rank 1 partial types. Then, either there exists a weakly-minimal formula that is almost
$\PP_0$-internal, or every complete type $p\in S(A)$ that is internal in $\PP_0$ is essentially
1-based over $\emptyset$ by means of $\tau^f$. In particular, either there exists a weakly-minimal
formula that is almost $\PP_0$-internal, or whenever $p\in S(A)$, where $A$ is countable, and
$p$ is analyzable in $\PP_0$, $p$ is analyzable in $\PP_0$ by
essentially 1-based types by means of $\tau^f$.
\end{theorem}

The most general dichotomy theorem of this type that we present is the following theorem that generalizes [S2, Theorem 2.3] to
any $\emptyset$-invariant family of $SU$-rank 1 types. The proof of this theorem is almost identical to the proof of [S2, Theorem 2.3]
(but note that the next version that we present with a proof contains all modifications that are needed for the proof of it).

\begin{theorem}\label{dichotomy thm}
Let $T$ be any hypersimple theory. Let $\Upsilon$ be a projection-closed family of topologies. Let
$\PP_0$ be an $\emptyset$-invariant family of $SU$-rank 1 types. Then, either there exists an unbounded
$\Upsilon$-open set (over some small set $A$) that is almost $\PP_0$-internal (and in particular has
finite $SU$-rank), or every complete type $p\in S(A)$ that is internal in $\PP_0$ is s-essentially
1-based over $\emptyset$ by means of $\Upsilon$. In particular, either there exists an unbounded
$\Upsilon$-open set that is almost $\PP_0$-internal, or whenever $p\in S(A)$ and $p$ is analyzable in $\PP_0$ ,
$p$ is analyzable in $\PP_0$ by s-essentially 1-based types
by means of $\Upsilon$.
\end{theorem}

The next theorem is a version of  Theorem \ref{dichotomy thm} for a countable language with a stronger consequence and is a generalization
of [S2, Theorem 2.11] to $\emptyset$-invariant family of $SU$-rank 1 types. We give the complete proof of this theorem.


\begin{theorem}\label{empty_inv_su1_dich_countable}
Let $T$ be any countable hypersimple theory. Let $\Upsilon$ be a projection-closed family of
topologies such that $\{ a\in \CC^x \vert a\not\in acl(A)\}\in \Upsilon_{x,A}$ for all $x$ and set
$A$ . Let $\PP_0$ be an $\emptyset$-invariant family of $SU$-rank 1 types. Then, either there exists an
unbounded type-definable $\Upsilon$-open set over some small set that is almost $\PP_0$-internal and
has \textbf{bounded} finite $SU$-rank, or every complete type $p\in S(A)$ that is internal in $\PP_0$
is essentially 1-based over $\emptyset$ by means of $\Upsilon$. In particular, either there exists
an unbounded type-definable $\Upsilon$-open set that is almost $\PP_0$-internal and has \textbf{bounded} finite
$SU$-rank, or whenever $p\in S(A)$, where $A$ is countable, and  $p$ is analyzable in $\PP_0$,
$p$ is analyzable in $\PP_0$ by essentially 1-based types by means of $\Upsilon$.
\end{theorem}

\noindent\proof $\Upsilon$ will be fixed and we'll freely
omit the phrase "by means of $\Upsilon$". To see the "In particular" part, work over a countable $A$ and assume
that every $p'\in S(A')$, with countable $A'\supseteq A$ , that is internal in $\PP_0$, is essentially 1-based
over $A$. Moreover, assume $p\in S(A)$ is non-algebraic and every non-algebraic extension of $p$ is
non-foreign to $\PP_0$. Then, for $a\models p$ there exists $a'\in dcl(Aa)\backslash acl(A)$ such
that $tp(a'/A)$ is $\PP_0$-internal and thus essentially 1-based over $A$ by our assumption. Thus,
by repeating this process we get that $p$ is analyzable in $\PP_0$ by essentially 1-based types.

We now prove the main part. Assume there exists $p\in S(A)$ that is internal in $\PP_0$, and $p$ is not
essentially 1-based over $\emptyset$. By the definition, there exist a finite tuple $d$ of
realizations of $p$ and $b$ that is independent from $d$ over $A$, and a finite tuple $\bar
c$ of realizations of types from $\PP_0$ over $Ab$ such that $d\in dcl(Ab\bar c)$, and there exists a type-definable $\Upsilon$-open set $\UU$ over
$Ad$ such that $a$ is independent from $A$ for all $a\in \UU$ and $\{a\in \UU \vert
Cb(a/Ad)\not\subseteq acl(a)\}$ is not nowhere dense in the Stone-topology of $\UU$. So, since
$\Upsilon$ refines the Stone-topology, by intersecting $\UU$ with a definable set, we may assume
that $\{a\in \UU \vert Cb(a/Ad)\not\subseteq acl(a)\}$ is dense in the Stone-topology of $\UU$.
\noindent Now, for each (finite) subsequence $\bar c_0$ of $\bar c$, let $$F_{\bar c_0}=\{ a\in \UU
\vert\ \exists b',\bar c'_0,\bar c'_1\ \mbox{s.t.}\ tp(b'\bar c'_0\bar c'_1/Ad)=tp(b\bar c_0(\bar
c\backslash \bar c_0)/Ad)\ \mbox{and} \nonforkempty{a}{Ab'\bar c'_0}\}.$$ Note that since $d$ is
independent from $b$ over $A$, any $a\in\UU$ is independent from $Ab'$ whenever
$tp(b'/Ad)=tp(b/Ad)$ and $\nonfork{a}{b'}{Ad}$. Thus $F_{\langle\rangle}=\UU$. Let $\bar c^*_0$ be
a maximal subsequence (with respect to inclusion) of $\bar c$ such that $F_{\bar c^*_0}$ has
non-empty Stone-interior in $\UU$ over $Ad$ (note that $F_{\bar c}$ has no Stone-interior
relatively in $\UU$). Let $\UU^*=\bigcap_{\bar c^*_0\subset\bar c'\subseteq\bar c} \UU\backslash
F_{\bar c'}$. Note that each $F_{\bar c'}$ is Stone closed relatively in $\UU$. Thus $\UU^*$ is
Stone-dense and Stone-open in $\UU$ and therefore there exists a
non-empty relatively Stone-open in $\UU$ set $W^*\subseteq F_{\bar c_0^*}\cap \UU^*$.
As $\UU$ is type-definable, we may  assume $W^*$ is type-definable.

\begin{claim}\label{subclaim0_main}
$W^*$ is a non-empty $\Upsilon$-open set over $Ad$ such that $\{a\in W^* \vert\
Cb(a/Ad)\not\subseteq acl(a)\}$ is dense in the Stone-topology of $W^*$ and for every $a\in W^*$ we
have: there exists $b'\bar c'_0\bar c'_1\models tp(b\bar c^*_0(\bar c\backslash \bar c^*_0)/Ad)$
such that $a$ is independent from $Ab'\bar c'_0$ over $\emptyset$ and moreover, for every $b'\bar
c'_0\bar c'_1\models tp(b\bar c^*_0(\bar c\backslash \bar c^*_0)/Ad)$ such that $a$ is independent
from $Ab'\bar c'_0$ we necessarily have $\bar c'_1\in acl(aAb'\bar c'_0)$.
\end{claim}

\proof The first part is immediate by the fact that $W^*\subseteq F_{\bar c_0^*}$. For the "moreover" part note that since $a\in W^*\subseteq \UU^*$, we get that by the definition of $\UU^*$,  $c'\in acl(aAb'\bar c'_0)$ for every $b'\bar c'_0\bar c'_1\models tp(b\bar c^*_0(\bar c\backslash \bar c^*_0)/Ad)$ and $c'\in \bar c'_1$ (as  $SU(c'/Ab')\leq 1$ for every $c'\in \bar c'_1$). \qed\\

\noindent Let us now define a set $V$ over $Ad$ by\\  $$V=\{(e',b',\bar c'_0,\bar c'_1,a') \vert \
\mbox{if}\ tp(b'\bar c'_0\bar c'_1/Ad)=tp(b\bar c^*_0(\bar c\backslash \bar c_0^*)/Ad)\ \mbox{and}
\nonforkempty{a'}{Ab'\bar c'_0}$$  $$\mbox{then}\ e'\in acl(Cb(Ab'\bar c'_0\bar c'_1/a'))\}.$$

\noindent Let $V^*=\{e' \vert \exists a'\in W^*\ \forall b',\bar c'_0,\bar c'_1\ V(e',b',\bar
c'_0,\bar c'_1,a')\}.$

\begin{claim}\label{subclaim1_main}
$V^*$ is a $\Upsilon$-open set over $Ad$.
\end{claim}

\proof Recall the following fact [S2, Proposition 2.4].

\begin{fact}\label {open Cb}
Let $q(x,y)\in S(\emptyset)$ and let $\chi(x,y,z)$ be an $\emptyset$-invariant set such that for
all $(c,b,a)\models \chi(x,y,z)$ we have $b\unrhd_a bc$. Then the set $$\UU=\{(e,c,b,a) \vert\ e\in
acl(Cb(cb/a))\}$$ is relatively Stone-open inside the set
$$F=\{(e,c,b,a)\vert\ \nonforkempty{b}{a}, \models\chi(c,b,a), tp(cb)=q\}.$$ (where
$e$ is taken from a fixed sort too).
\end{fact}

\noindent By Fact \ref{open Cb} and Claim \ref{subclaim0_main}, there exists a Stone-open set $V'$
over $Ad$ such that for all $a'\in W^*$ and for all $e',b',\bar c'_0,\bar c'_1$ we have
$V'(e',b',\bar c'_0,\bar c'_1,a')$ if and only if $V(e',b',\bar c'_0,\bar c'_1,a')$. Thus, we
may replace $V$ by $V'$ in the definition of $V^*$. As Stone-open sets are closed under the
$\forall$ quantifier, the $\Upsilon$ topology refines the Stone-topology and closed under product
by a full Stone-space and closed under projections, we conclude that $V^*$ is a $\Upsilon$-open
set.$\ \ \ \ \Box$\\

\begin{claim}\label{subclaim2_main}
For appropriate sort for $e'$, the set $V^*$ is unbounded and is almost $\PP_0$-internal (over $Ad$)
and thus has finite $SU$-rank over $Ad$.
\end{claim}

\proof First, note the following general observation.

\begin{fact}\label{dcl_cb remark}
Assume $d\in dcl(c)$. Then $Cb(d/a)\in dcl(Cb(c/a))$ for all $a$.
\end{fact}

\noindent Let $a^*\in W^*$ be such that $Cb(a^*/Ad)\not\subseteq acl(a^*)$. Then
$Cb(Ad/a^*)\not\subseteq acl(Ad)$. By Fact \ref{dcl_cb remark}, there exists $e^*\not\in acl(Ad)$
such that $e^*\in acl(Cb(Ab'\bar c'_0\bar c'_1/a^*))$ for all $b'\bar c'_0\bar c'_1\models tp(b\bar
c^*_0(\bar c\backslash \bar c^*_0)/Ad)$. In particular, $e^*\in V^*$. Thus, if we fix the sort for
$e'$ in the definition of $V^*$ to be the sort of $e^*$, then $V^*$ is unbounded. Now, let $e'\in
V^*$. Then for some $a'\in W^*$, $\models V(e',\bar c'_0,\bar c'_1,b',a')$ for all $b',\bar
c'_0,\bar c'_1$. By Claim \ref{subclaim0_main}, there exists $b'\bar c'_0\bar c'_1\models
tp(b\bar c^*_0(\bar c\backslash \bar c^*_0)/Ad)$ such that $a'\ \mbox{is\ independent\ from}\
Ab'\bar c'_0\ \mbox{over}\ \emptyset$. Thus, by the definition of $V^*$ and $V$, $e'\in
acl(Cb(Ab'\bar c'_0\bar c'_1/a'))$. Since $Ab'$ is independent from $a'$ over $\emptyset$, $tp(e')$
is almost $\PP_0$-internal (as $Cb(Ab'\bar c'_0\bar c'_1/a')$ is in the definable closure of any
Morley sequence of $Lstp(Ab'\bar c'_0\bar c'_1/a')$ ). In particular, $tp(e'/Ad)$ is almost
$\PP_0$-internal by Fact \ref{some internal} and therefore $tp(e'/Ad)$ has finite $SU$-rank. $\ \ \ \ \Box$\\

\begin{claim}
There exists $V^{**}\subseteq V^*$ that is unbounded, type-definable and $\Upsilon$-open over $Ad$.
\end{claim}

\proof By the definition of $V^*$ and the proof of
Claim \ref{subclaim1_main} there exist a Stone open set $V_0$ over $Ad$ such that $V^*=\{ e'
\vert \exists a'\in W^*\ (V_0(e',a'))\}$.  By replacing $V_0$ by a definable set and using the fact that $W^*$ is type-definable and that $\Upsilon$ is a
projection-closed family of topologies we get the required set $V^{**}$ $\ \ \ \ \Box$\\

\noindent Now, by the proof of Claim \ref{subclaim2_main} we know that for all $e'\in V^{**}$ we have
$e'\in acl(Cb(Ab'\bar c'_0\bar c'_1/a'))$ for some $a'\in W^*$ and some $b',\bar c'_0,\bar c'_1$
such that $a'$ is independent from $Ab'\bar c'_0\ \mbox{over}\ \emptyset$ and $b'\bar c'_0\bar
c'_1\models tp(b\bar c^*_0(\bar c\backslash \bar c^*_0)/Ad)$. Let $q=tp(Ab)$. For every
$\chi=\chi(x,y_0,...,y_n,\bar z_0,\bar z_1,...\bar z_n)\in L$ (for some $n<\omega$) such that $\forall y_0 y_1 ...y_n\bar z_0\bar z_1...\bar z_n  \ \exists^{<\infty}x\ \chi(x,y_0,y_1,...y_n,\bar z_0,\bar z_1,...\bar z_n)$, let

$$F_{\chi}=\{e\in V^{**} \vert\ \models\chi(e,C_0,C_1,..C_n,\bar c_0,\bar c_1,...\bar c_n)\  \mbox{for\ some\ } \bar c_0,...\bar c_n\ \mbox{and\ some}$$
$$\ \ \ \ \ \ \ \ \ \ \ \ \ \ \ \emptyset-\mbox{independent\ sequence}  (C_i\vert i\leq n)\ \mbox{of\ realization of}\ q$$  $$ \ \ \ \ \ \ \ \ \ \ \ \ \ \ \ \mbox{such that}\  tp(C_i,\bar c_i)=tp(Ab,\bar c)\ \mbox{and}\ \nonforkempty{e}{(C_i\vert i\leq n)} \}.$$

\noindent Note that eaxh $F_\chi$ is type-definable. By the aforementioned, we get that $V^{**}\subseteq \bigcup_{\chi} F_{\chi}$ (the
union is over each $\chi$ as above). By the Baire category theorem applied to the Stone-topology
of the Stone-closed set $V^{**}\backslash acl(Ad)$, there exists $\theta\in L(Ad)$ such that
$$\tilde V\equiv \theta^\CC\cap(V^{**}\backslash acl(Ad))\neq\emptyset\ \mbox{and\ } \tilde
V\subseteq F_{\chi^*}$$ for some $\chi^*$ as above. Clearly, $\tilde V$ is unbounded,
type-definable and $\Upsilon$-open (by the assumptions on $\Upsilon$). Now, there exists a fixed $m^*<\omega$ such that for every $a\in \tilde
V$, $SU(a/Ad)\leq m^*$ and $tp(a/Ad)$ is almost $\PP_0$-internal (as $tp(a)$ is almost
$\PP_0$-internal). This completes the proof of the main part of the theorem. $\ \ \ \ \Box$\\

The proof of the main result of this section now follows exactly in the same way as in [S2].
We write the proof for completeness.\\

\noindent\textbf{Proof of Theorem \ref{empty_inv_su1_dich_pcft}}
Our assumptions are clearly a special case of the assumptions of Theorem \ref{empty_inv_su1_dich_countable}, thus
we only need to prove the first part. By the conclusion of Theorem \ref{empty_inv_su1_dich_countable}, we may assume that there
exists a $\tau^f$-open set $\UU$ of bounded finite $SU$-rank over some small set $A$ that is almost
$\PP_0$-internal. By Fact \ref{tau bounded SU}, there exists exists a weakly-minimal $\theta(x,b)\in L(B)$ for some
small set $B\supseteq A$, such that $\theta^\CC\subseteq \UU\cup acl(B)$. Now, $tp(a/B)$ is almost
$\PP_0$-internal for every $a\in \theta^\CC$, and so $tp(a/b)$ ($b$ is the parameter of $\theta(x,b)$
) is almost $\PP_0$-internal for every $a\in \theta^\CC$ (by taking non-forking
extensions).$\ \ \ \ \Box$\\

\section{Main result}
From now on we assume $T$ is an arbitrary simple theory with elimination of imaginaries unless
stated otherwise. We work in $\CC$.

\begin{definition}\label {Coord_finite_rank_def}\em
1) We say that $T$ is \em analyzable in $SU$-rank 1 types \em if every type is analyzable in the family of $SU$-rank 1 types.

\noindent 2) We say that $T$ is \em coordinatized in finite rank \em if for every $a\in\CC$ and $A\subseteq
\CC$ such that $a\not\in acl(A)$ there exists $a'\in acl(aA)\backslash acl(A)$ with
$SU(a'/A)<\omega$.
\end{definition}

\begin{lemma}\label{cifr}
Assume $T$ is hypersimple. $T$ is coordinatized in finite rank iff $T$ is analyzable in $SU$-rank 1 types.
\end{lemma}

\proof If $T$ is analyzable in $SU$-rank 1 types then clearly $T$ is
coordinatized in finite rank. Assume now that $T$ is coordinatized in finite rank. We first note
the following.

\begin{claim}\label{claim_non_orth_su1}
Let $T$ be any simple theory. Let $a\in\CC$ be such that $SU(a)=n<\omega$ and such that for some
$b\in \CC$ with $SU(b)<\infty$ we have $SU(a/b)=n-1$. Then $tp(a)$ is non-orthogonal to an $SU$-rank 1
hyperimaginary type.
\end{claim}

\proof Let $e=Cb(Lstp(a/b))$ ($e$ is a hyperimaginary). Since $SU(e)<\infty$ (as we assume
$SU(b)<\infty$), there exists a set $A$ such that $SU(e/A)=1$. By extension we may clearly assume
$\nonfork{a}{A}{e}$. We claim that $e\in bdd(aA)$  (*). Indeed, otherwise $\nonfork{e}{a}{A}$ and so $e\in bdd(A)$
 (as $tp(a/e)$ is canonical), a contradiction to $SU(e/A)=1$.  Now, $SU(a/eA)=SU(a/e)=n-1$. By (*), $SU(a/A)=SU(ae/A)\geq SU(a/eA)+SU(e/A)=n$.
 Thus $\nonforkempty{a}{A}$ and so $tp(a)$ is non-orthogonal to $tp(e/A)$. \qed\\

\noindent Now, let $a,A$ be given such that $a\not\in acl(A)$.  By our assumption, there exists $a'\in acl(aA)\backslash acl(A)$ with
$SU(a'/A)=n$ for some $n<\omega$. Let $b\in\CC$ be such that $SU(a/Ab)=n-1$ and let $(b_i\vert i<\alpha)$ be such that $b_i\in acl(bA)$
 and\\ $0<SU(b_i/Ab_{<i})<\omega$ for all $i<\alpha$ and such that $acl(Ab)=acl(A\cup\{b_i\vert i<\alpha\})$.  As $\fork{a'}{b}{A}$,
 there exists a minimal $i^*<\alpha$ such that  $\fork{a'}{ \{b_i\vert  i\leq i^* \} }{A}$. Then  $\nonfork{a'}{ \{b_i\vert  i< i^* \} }{A}$.
 Now, $a'$ satisfies the assumptions of  Claim \ref{claim_non_orth_su1} when working over  $A\cup \{b_i\vert  i< i^* \}$.
 Thus $tp(a'/A\cup \{b_i\vert  i< i^* \})$ is non-orthogonal to an $SU$-rank 1 type which we may clearly assume to be a type of an imaginary.
 Thus $tp(a'/A)$ is non-orthogonal an $SU$-rank 1 imaginary type and so is $tp(a/A)$ (as $a'\in acl(aA)$).\\

We start with the following proposition that generalizes the main result in [S1]; the proof is
similar but applies Theorem \ref{empty_inv_su1_dich_pcft}.

\begin{proposition}\label{reduction_prop}
Assume $T$ is a countable hypersimple theory that is coordinatized in finite rank
and eliminates $\exists^\infty$. Then there exists a weakly minimal formula.
\end{proposition}

\proof By Lemma \ref{cifr}, $T$ is analyzable in SU-rank 1 types. As $T$ eliminates
$\exists^\infty$, $T$ is EPFO by Fact \ref {fact 1}. By Fact \ref{ext pcft}, $T$ is PCFT. Let $\PP_0$ be the
family of all SU-rank 1 types. By Theorem \ref{empty_inv_su1_dich_pcft}, we may assume that every complete finitary
type over a countable set is analyzable in $\PP_0$ by essentially 1-based types by means of
$\tau^f$. We recall the following fact that for simplicity we state for a special case [S3,
Corollary 4.5]. Infact, the result is valid for a large class of sets ($\tilde\tau^f_{low}$-sets)
instead of the specific set $\UU_0$ that is fixed in the following statement.

\begin{fact}\label{cor1}
Let $T$ be a countable simple theory with EPFO. Let $\UU_0=\CC^s\backslash acl(\emptyset)$ for some
non-algebraic sort $s$ say. Assume for every $a\in\UU_0$ there exists $a'\in acl(a)\backslash
acl(\emptyset)$ such that $SU_{se}(a')<\omega$. Then there exists an unbounded
$\tau_{\infty}^f$-open set $\UU$ over a finite set such that $\UU$ has bounded finite
$SU_{se}$-rank.
\end{fact}

\noindent $T$ satisfies the assumptions of Fact \ref{cor1}, so let $\UU$ be a set as in its conclusion. In particular, $SU_s(\UU)=n$ for some
$n<\omega$. Recall now the following easy lemma.

\begin{lemma}\label{su_s finite to 1} \em [S1, Lemma 7.4] \em
Assume $\UU$ is an unbounded $\tau^f_\infty$-open set of bounded finite $SU_s$-rank over some
finite set $A$. Then there exists a $\tau^f_\infty$-open set $\UU^*\subseteq \UU$ over some finite
set $B^*\supseteq A$ of $SU_s$-rank 1.
\end{lemma}

\noindent By Lemma \ref{su_s finite to 1}, we may assume $SU_s(\UU)=1$, $\UU$ is a type-definable
$\tau^f$-open set over a finite set $A_0$. We claim $SU(\UU)=1$. Indeed, otherwise there
exists $a$ and $d\in \UU$ such that $\fork{d}{a}{A_0}$ and $d\not\in acl(aA_0)$. Since every
finitary type over a countable set is analyzable in $\PP_0$, there exists $(a_i \vert
i\leq\alpha)\subseteq dcl(aA_0)$ with $a_\alpha=a$ (where $\alpha<\omega_1$) such that
$tp(a_i/A_0\cup\{a_j \vert j<i\})$ is essentially 1-based over $A_0$ by means of $\tau^f$ for all
$i\leq\alpha$. Now, let $i^*\leq\alpha$ be minimal such that there exists $d'\in\UU$ satisfying
$\fork{d'}{\{a_i \vert i\leq i^*\}}{A_0}$ and $d'\not\in acl(A_0\cup \{a_i \vert i\leq i^*\})$.
Pick $\phi(x,a')\in L(A_0\cup \{a_i \vert i\leq i^*\})$ that forks over $A_0$ and such that
$\phi(d',a')$. Let
$$V=\{d\in\UU \vert\ \phi(d,a')\ \mbox{and}\ d\not\in acl(A_0\cup \{a_i \vert i\leq i^*\})\ \}.$$
By minimality of $i^*$, $d$ is independent from $\{a_i \vert i<i^*\}$ over $A_0$ for all $d\in V$.
Clearly $V$ is type-definable and by Fact \ref{tau extensions}, $V$ is a $\tau^f$-open set over
$A_0\cup \{a_i \vert i\leq i^*\}$. Now, since $tp(a_{i^*}/A_0\cup\{a_i \vert i<i^*\})$ is
essentially 1-based over $A_0$ by means of $\tau^f$, the set
$$\{d\in V \vert\ Cb(d/A_0\cup \{a_i \vert i\leq i^*\})\in bdd(dA_0)\}$$ contains a relatively
Stone-open and Stone-dense subset of $V$. In particular, there exists $d^*\in V$ such that
$tp({d^*}/{A_0\cup \{a_i \vert i\leq i^*\}}$ doesn't fork over ${acl(A_0d^*)\cap acl(A_0\cup\{a_i
\vert i\leq i^*\})}$. Since we know $\fork{d^*}{A_0\cup\{a_i \vert i\leq i^*\}}{A_0}$, Fact
\ref{stable forking} implies $\sfork{d^*}{A_0\cup\{a_i \vert i\leq i^*\}}{A_0}$. Hence $d^*\in V$
implies $SU_s(d^*/A_0)\geq 2$, which contradicts $SU_s(\UU)=1$. Thus we have proved $SU(\UU)=1$.
Now, by Fact \ref{tau bounded SU} there exists a definable set of $SU$-rank 1.\qed\\

Before stating the main theorem, we give some terminology and easy remarks. Recall that we work in $\CC=\CC^{eq}$ and that
$(\RR^\lambda_{\CC},\leq_{\CC})$ is the partial order of reducts of $\CC$ of size$\le\lambda$.

\begin{definition}\em
Let $\CC\vert L^-\in \RR^\lambda_{\CC}$. We will say that $\CC\vert L^-$ is \em eq-closed \em if $T^-=Th(\CC\vert L^-)$ has uniform elimination of imaginaries, i.e.
for every definable equivalence relation $E\in L^-$ on $S_0\times
S_1\times...S_k$, where $S_i$ are sorts of $L^-$, there is a definable function $f_E\in L^-$ whose
domain is $(S_0\times S_1\times...S_k)^\CC$ such that for all $\bar a,\bar b$, we have  $f_E(\bar a)=f_E(\bar b)$ iff $E(\bar a, \bar b)$.
\end{definition}

\begin{remark}\label{eq_club}\em
\noindent  For every reduct  $\CC\vert L^-\in \RR^\lambda_{\CC}$ there exists a reduct $\CC\vert L^*\in\RR^\lambda_{\CC} $ that is eq-closed and is an expansion of $\CC\vert L^-$.
 Thus for every infinite cardinal $\lambda$,  the set of reducts in $\RR^\lambda_{\CC}$ that are eq-closed is a club in $(\RR^\lambda_{\CC},\leq_{\CC})$.
\end{remark}

\proof Expand the reduct $\CC\vert L^-$ of $\CC=\CC^{eq}$ by adding for every definable equivalence relation $E$ on $S_0\times
S_1\times...S_k$, where $S_i$ are sorts of $L^-$ and $E\in L^-$, a definable function $f_E\in L$ whose
domain is $(S_0\times S_1\times...S_k)^\CC$ and is onto the interpretation of some sort of $L$ such that $f_E(\bar a)=f_E(\bar b)$ iff $E(\bar a, \bar b)$.
Now, the resulting expansion will have uniform elimination of imaginaries. It is immediate that the set of eq-closed reducts in $\RR^\lambda_{\CC}$  is closed in $(\RR^\lambda_{\CC},\leq_{\CC})$.\\

\noindent Now, note the following easy general remark on clubs.

\begin{remark}\label{club_intersection}\em
Let $(P,\leq)$ be a directed partial order that is $\omega$-closed (i.e. any increasing sequence $(a_i\vert i<\omega)$ has a supremum).
Then the intersection of finitely many clubs in $(P,\leq)$ is a club.
\end{remark}

In the proof we will refer to the following notion.

\begin{definition}\em
We say that $T$ is strongly non-supersimple if $D(\phi(x,a))=\infty$ for every non-algebraic
$\phi(x,a)\in L(\CC)$.
\end{definition}

\begin{remark}\em
 Note that  $T$ is strongly non-supersimple iff  for every non-algebraic $\phi(x,a)\in L(\CC)$ there exists a non-algebraic $\psi(x,b)\in L(\CC)$ such that
 $\psi(x,b)\vdash\phi(x,a)$ and $\psi(x,b)$ forks over $a$ iff there does not exist a weakly minimal formula.
\end{remark}

\begin{theorem}\label{main_thm}
Let $T=T^{eq}$ be a hypersimple unidimensional theory. Assume there is a club of countable reducts
of $T$ in $(\RR^{\aleph_0}_{\CC},\leq_{\CC})$ that are coordinatized in finite rank. Then $T$ is supersimple.
\end{theorem}

\proof First, if $T$ is not strongly non-supersimple then we are done by Fact \ref{supersimple definable}. Therefore, we may assume $T$
is strongly non-supersimple. By Fact \ref{elimination_exists_infty}, $T$ eliminates $\exists^\infty$, thus every
reduct of $T$ eliminates $\exists^\infty$.


\begin{claim}\em\label{main_claim}
The set $\tilde \CC_1$ of countable strongly non-supersimple reducts of $\CC$ is a club in $(\RR^{\aleph_0}_{\CC},\leq_{\CC})$.
\end{claim}

\proof First, we prove that $\tilde \CC_1$ is unbounded in $(\RR^{\aleph_0}_{\CC},\leq_{\CC})$.
 Let $\CC\vert L^-\in \RR^{\aleph_0}_{\CC}$.  We construct by induction an increasing sequence of reducts
$(\CC_n\vert n<\omega)$, $\CC_n\in \RR^{\aleph_0}$, where $\CC_n=\CC\vert L^-_n$ for some countable sublanguage $L^-_n$ of $L$, $T^-_n=Th(\CC_n)$ in the following way.
Let $\CC_0=\CC\vert L^-_0$, $L^-_0=L^-$ and assume $L^-_k$ have already been defined for $k\leq n$. We define $L^-_{n+1}$. For any fixed $\phi(x,y)\in L^-_n\backslash L^-_{n-1}$
we define a finite set of formulas $\Delta_\phi=\{\psi_i \vert i\leq n(\phi)\}$, where
$\psi_i=\psi_i(x,y_i)\in L= L(T)$, $n(\phi)<\omega$ in the following way. Since $T$ is strongly
non-supersimple, for every $a\in\CC$ such that $\exists^\infty x\phi(x,a)$, there exists
$\psi_a(x,z)\in L$ and some $b\in\CC$ such that $\psi_a(x,b)\vdash \phi(x,a)$, $\psi_a(x,b)$ forks
over $a$ and $\exists^\infty x\psi_a(x,b)$ (*). For every $\psi(x,z)\in L$ let
$$\theta_{\phi,\psi}(z,y)=\exists^\infty x \psi(x,z)\wedge \forall
x(\psi(x,z)\rightarrow\phi(x,y)).$$

\noindent By elimination of $\exists^\infty$ (in $\CC$), $\theta_{\phi,\psi}(z,y)$ is definable.
Now, let $Q_{\psi,\theta_{\phi,\psi}}(y)$ be the relation in Fact \ref{foext} defined for
$\theta_{\phi,\psi},\psi$. So, for every $a\in\CC$, $\neg Q_{\psi,\theta_{\phi,\psi}}(a)$ iff there
exists $b$ such that $\psi(x,b)$ is not algebraic, $\psi(x,b)\vdash \phi(x,a)$ and $\psi(x,b)$
forks over $a$. Since $T$ is EPFO, we know that each $\neg Q_{\psi,\theta_{\phi,\psi}}$ is
Stone-open. By (*), in $\CC$: $$\exists^\infty x\phi(x,y)\vdash \bigvee_{\psi\in L}\neg
Q_{\psi,\theta_{\phi,\psi}}(y).$$ By compactness, there are $\psi_0=\psi_0(\phi), ...,
\psi_{n(\phi)}=\psi_{n(\phi)}(\phi)\in L$ such that in $\CC$: $$\exists^\infty x\phi(x,y)\vdash
\bigvee_{i<n(\phi)}\neg Q_{\psi_i,\theta_{\phi,\psi_i}}(y)\ \ (**).$$   Now, let $\Delta_\phi=\{\psi_i(\phi)
\vert\ i\leq n(\phi)\}$. Let $L^-_{n+1}$  be the set of formulas generated by the set
$$\nu_{n+1}=L^-_n\cup\bigcup\{\Delta_\phi\vert \phi(x,y)\in L^-_n\backslash L^-_{n-1}\},$$ where the set of
sorts attached to $L^-_{n+1}$ is the set of all sorts of variables that appears in $\nu_{n+1}$. Let $\CC_{n+1}=\CC\vert L^-_{n+1}$.
Now, let $L^-_\omega=\bigcup_{n<\omega} L^-_n$ and let $\CC_\omega=\CC\vert L^-_\omega$, $T^-_\omega=Th(\CC_\omega)$.
We claim that $T^-_\omega$ is strongly non-supersimple. Indeed, given a formula $\phi_\omega(x,y)\in L^-_\omega$,  let $a\in \CC_\omega$ be
such that $\models \exists^\infty x \phi_\omega(x,a)$. Then, by (**) there exists $\psi(x,z)\in \Delta_{\phi_\omega}\subseteq L^-_\omega$ such that for some $b$
we have $\psi(x,b)\vdash \phi_\omega(x,a)$ and $\psi(x,b)$ forks over $a$ in $\CC$ and thus in particular $\psi(x,b)$ forks over $a$ in $\CC_\omega$.
Thus $T^-_\omega$ is strongly non-supersimple. Now, to show that $\tilde\CC_1$  is closed  in $(\RR^{\aleph_0}_{\CC},\leq_{\CC})$,
let $\tilde\CC$ be a chain in $\tilde\CC_1$. We claim that $\CC^*=\bigcup\tilde \CC$
(the universe of $\CC^*$ is the union of the interpretations of the sorts of all members of $\tilde\CC$ and likewise for the definable sets of $\CC^*$)  is strongly non-supersimple.
 Indeed, let $\phi(x,a)\in L(\CC^*)$ be non-algebraic. Then there exists $\CC_0=\CC\vert L_0\in \tilde\CC$ for some countable sublanguage $L_0$ of $L$ such that $\phi(x,a)\in L_0(\CC_0)$.
 Since $Th(\CC_0)$ is strongly non-supersimple, there exists a non-algebraic $\psi(x,b)\in L_0(\CC_0)$ such that $\psi(x,b)\vdash \phi(x,a)$ and $\psi(x,b)$ forks over $a$ in $\CC_0$.
 By Ramsey and compactness, there exists a formula $\psi(x,b')$  that is $a$-conjuagate to $\psi(x,b)$ in $\CC_0$ and that forks over $a$ in the sense of $\CC^*$.
 Thus $Th(\CC^*)$ is strongly non-supersimple.\qed\\

 \noindent By Claim \ref{main_claim},  Claim \ref{reduct_elim},  Remark \ref{eq_club}, the assumptions of the theorem and Remark \ref{club_intersection}, there is a club of reducts in $\RR^{\aleph_0}_{\CC}$ that are strongly non-supersimple,
 hypersimple, eq-closed and coordinatized in finite rank. Any such reduct contradicts Proposition \ref{reduction_prop} (as elimination of $\exists^\infty$ is clearly preserved under reducts).\qed

Ziv Shami, E-mail address: zivsh@ariel.ac.il.

\end{document}